\documentclass[11pt]{amsart}
\usepackage{amsmath}
\usepackage{amsfonts}
\usepackage{amsthm}
\usepackage{amssymb}
\usepackage{tikz-cd}
\usepackage{graphics}
\usepackage{array}
\usepackage{enumerate}
\usepackage{dsfont}
\usepackage{color}
\usepackage{wasysym}
\usepackage{hyperref}
\usepackage{pdfsync}
\usepackage{mathrsfs}

\newcommand{\N}{{\mathbb N}}

\newcommand{\C}{{\mathbb C}}
\newcommand{\Q}{{\mathbb Q}}

\DeclareMathOperator{\Alb}{Alb}

\DeclareMathOperator{\Pic}{Pic}

\DeclareMathOperator{\codim}{codim}
\DeclareMathOperator{\Supp}{Supp}

\DeclareMathOperator{\IT}{IT_0}
\DeclareMathOperator{\Exc}{Exc}

\DeclareMathOperator{\red}{red}
\DeclareMathOperator{\Var}{Var}
\DeclareMathOperator{\rank}{rank}

\theoremstyle{plain}

\newtheorem{thm}{Theorem}[section]
\newtheorem{conj}[thm]{Conjecture}
\newtheorem{prop}[thm]{Proposition}
\newtheorem{coro}[thm]{Corollary}
\newtheorem{lemma}[thm]{Lemma}

\newtheorem*{Iitaka conjecture}{Iitaka conjecture}
\newtheorem*{Viehweg conjecture}{Viehweg conjecture}

\theoremstyle{definition}

\newtheorem{defi}[thm]{Definition}

\newtheorem*{ac}{Acknowledgements}

\begin{document}

\title[On surjective morphisms to abelian varieties]{On surjective morphisms to abelian varieties and a generalization of the Iitaka conjecture}
	
	\author{Fanjun Meng}
	\address{Department of Mathematics, Johns Hopkins University, 3400 N. Charles Street, Baltimore, MD 21218, USA}
	\email{fmeng3@jhu.edu}


\begin{abstract}
We explore the relationship between fibrations arising naturally from a surjective morphism to an abelian variety. These fibrations encode geometric information about the morphism. Our study focuses on the interplay of these fibrations and presents several applications. Then we propose a generalization of the Iitaka conjecture which predicts an equality of Kodaira dimension of fibrations, and prove it when the base is a smooth projective variety of maximal Albanese dimension.
\end{abstract}

	\maketitle
	\tableofcontents

\section{Introduction}\label{1}

In this paper, we explore the relationship between fibrations arising naturally from a surjective morphism to an abelian variety over $\C$. These fibrations encode geometric information about the morphism. Our study focuses on the interplay of these fibrations and presents several applications. To enhance clarity, we begin by discussing the applications before introducing the main theorems.

\subsection{Applications}

By applying the techniques of this paper, we give an equality of Kodaira dimension when the base is an abelian variety.

\begin{thm}\label{eq1}
Let $f\colon X \to A$ be a surjective morphism from a klt pair $(X, \Delta)$ to an abelian variety $A$, $F$ the general fiber of $f$, $m>1$ a rational number, and $D$ a Cartier divisor on $X$ such that $D\sim_{\Q}m(K_X+\Delta)$ and $f_*\mathcal{O}_X(D)\neq0$. For every positive integer $l$, we have
\begin{align*}
\kappa(X, K_X+\Delta)&=\kappa(F, K_F+\Delta|_F)+\kappa(A, \widehat{\det} f_*\mathcal{O}_X(lD))\\
&=\kappa(F, K_F+\Delta|_F)+\dim V^0(A, f_*\mathcal{O}_X(lD)).
\end{align*}
\end{thm}

The line bundle $\widehat{\det} f_*\mathcal{O}_X(lD)$ is the reflexive hull of $\det f_*\mathcal{O}_X(lD)$. Let $\mathcal{F}$ be a coherent sheaf on an abelian variety $A$. The cohomological support locus $V^0(A, \mathcal{F})$ is defined by
$$V^0(A, \mathcal{F})=\{\alpha\in\Pic^0(A)\mid\dim H^0(A, \mathcal{F}\otimes\alpha)>0\},$$
see also Definition \ref{csl}.

It is known that $\kappa(X, K_X+\Delta)\ge\kappa(F, K_F+\Delta|_F)$ by \cite{CP17}, which is the subadditivity of Kodaira dimension of fibrations over abelian varieties. It is also known that we have the stronger inequality $\kappa(X, K_X+\Delta)\ge\kappa(F, K_F+\Delta|_F)+\dim V^0(A, f_*\mathcal{O}_X(lD))$ if $m\ge1$ by \cite{Men22}. Theorem \ref{eq1} finds the difference between $\kappa(X, K_X+\Delta)$ and $\kappa(F, K_F+\Delta|_F)$, and shows that the inequalities in \cite[Theorem 1.6]{Men22} are actually equalities if $m>1$. When $X$ is smooth and $\Delta=0$, an immediate corollary is the following.

\begin{coro}\label{eq1c}
Let $f\colon X \to A$ be a surjective morphism from a smooth projective variety $X$ to an abelian variety $A$ and $F$ the general fiber of $f$. For every integer $m>1$ such that $f_* \omega_{X}^{\otimes m}\neq0$, we have
$$\kappa(X)=\kappa(F)+\kappa(A, \widehat{\det} f_* \omega_{X}^{\otimes m})=\kappa(F)+\dim V^0(A, f_* \omega_{X}^{\otimes m}).$$
\end{coro}

By estimating the dimension of $V^0(A, f_*\mathcal{O}_X(lD))$, we have the following corollary of Theorem \ref{eq1}.

\begin{coro}\label{eq1c2}
Let $g\colon X\to Y$ be a smooth model of the Iitaka fibration associated to $K_X+\Delta$ where $(X, \Delta)$ is a klt pair and $Y$ is a smooth projective variety. Let $f\colon X \to A$ be a surjective morphism to an abelian variety $A$ with general fiber $F$. Then 
$$\kappa(X, K_X+\Delta)\le\kappa(F, K_F+\Delta|_F)+q(Y).$$
\end{coro}

Given a projective variety $Y$, $q(Y)$ denotes the irregularity of $Y$, see Definition \ref{irregularity}. It is known that
$$\kappa(X, K_X+\Delta)\ge\kappa(F, K_F+\Delta|_F)+\dim A-q(G)$$
where $G$ is the general fiber of $g$ by \cite[Corollary 1.7]{Men22}. Corollary \ref{eq1c2} gives an upper bound of $\kappa(X, K_X+\Delta)$. In particular, if $q(Y)=0$, we have
$$\kappa(X, K_X+\Delta)=\kappa(F, K_F+\Delta|_F).$$

We have a quick corollary of Theorem \ref{eq1} if the Albanese morphism of $X$ is a surjective morphism. In Corollary \ref{eq1c2}, we let $f$ be the Albanese morphism of $X$. Then we deduce that
$$\kappa(X, K_X+\Delta)=\kappa(F, K_F+\Delta|_F)+q(Y)$$
by Theorem \ref{eq1} and \cite[Theorem 1.5]{Men21}.

With more work, we can generalize Theorem \ref{eq1} slightly.

\begin{thm}\label{eq2}
Let $f\colon X \to A$ be a surjective morphism from a klt pair $(X, \Delta)$ to an abelian variety $A$, $F$ the general fiber of $f$, $m>1$ a rational number, and $D$ a Cartier divisor on $X$ such that $D\sim_{\Q}m(K_X+\Delta)$ and $f_*\mathcal{O}_X(D)\neq0$. Let $L$ be a divisor on $A$ such that $\kappa(A, L)\ge0$. For every positive integer $l$, we have
\begin{align*}
\kappa(X, K_X+\Delta+f^*L)&=\kappa(F, K_F+\Delta|_F)+\kappa(A, \widehat{\det} f_*\mathcal{O}_X(lD)\otimes\mathcal{O}_A(L))\\
&=\kappa(F, K_F+\Delta|_F)+\dim V^0(A, f_*\mathcal{O}_X(lD)\otimes\mathcal{O}_A(L)).
\end{align*}
\end{thm}

When $X$ is smooth, $\Delta=0$, and $f$ is a fibration, Kawamata proved that
$$\kappa(X, K_X+f^*L)\ge\kappa(F, K_F)+\max\{\kappa(A, L),\,\Var(f)\}$$
in \cite[Theorem 1.1(ii)]{Kaw85a} under the assumption that the general fiber $F$ has a good minimal model. For the definition of the variation $\Var(f)$, see \cite[Section 1]{Kaw85a}. Theorem \ref{eq2} recovers his result when the base is an abelian variety since
$$\kappa(A, \widehat{\det} f_* \omega_{X}^{\otimes m})\ge\Var(f)\ge0$$
for some integer $m>1$ by \cite[Theorem 1.1(i)]{Kaw85a} and thus
$$\kappa(A, \widehat{\det} f_* \omega_{X}^{\otimes m}\otimes\mathcal{O}_A(L))\ge\max\{\kappa(A, L),\,\Var(f)\}.$$

In a different but related direction, we have the following corollary as an application of Theorems \ref{m1} and \ref{m2}.

\begin{coro}\label{ineq1}
Let $f\colon X \to A$ be a surjective morphism from a smooth projective variety $X$ to an abelian variety $A$ where $f$ is smooth over an open set $V\subseteq A$, and $m$ a positive integer. Then
$$\kappa(V)\ge\kappa(A, \widehat{\det} f_* \omega_{X}^{\otimes m})\ge \dim V^0(A, f_*\omega_X^{\otimes m}).$$
If $m>1$, then $\kappa(A, \widehat{\det} f_* \omega_{X}^{\otimes m})=\dim V^0(A, f_*\omega_X^{\otimes m})$.
\end{coro}

Given a smooth quasi-projective variety $V$, $\kappa(V)$ denotes the log Kodaira dimension, defined as follows: for any smooth projective compactification $Y$ of $V$ such that $D = Y\setminus V$ is a divisor with simple normal crossing support, we have $\kappa (V) = \kappa (Y, K_Y+D)$. It is independent of the choice of the compactifications of $V$, see \cite[Section 3]{Iit77}.

If $\kappa(V)=0$, Corollary \ref{ineq1} is known by \cite[Theorem B]{MP23}. By \cite[Theorem 3.6, Remark 3.7 and Lemma 3.8]{Men22}, it is already known that (even for klt pairs) $\kappa(A, \widehat{\det} f_* \omega_{X}^{\otimes m})\ge \dim V^0(A, f_*\omega_X^{\otimes m})$, and $\kappa(A, \widehat{\det} f_* \omega_{X}^{\otimes m})=\dim V^0(A, f_*\omega_X^{\otimes m})$ if $m>1$ under the hypotheses of Corollary \ref{ineq1}. If the morphism $f$ is a fibration, we also know that $\kappa(V)\ge\kappa(A, \widehat{\det} f_* \omega_{X}^{\otimes m})$ by \cite[Theorem 1.1]{Men22}. Corollary \ref{ineq1} shows that $\kappa(V)\ge\kappa(A, \widehat{\det} f_* \omega_{X}^{\otimes m})$ when $f$ is a surjective morphism but not necessarily a fibration, which could not be obtained by the method in \cite{Men22}. Instead, it is obtained as a corollary of Theorems \ref{m1} and \ref{m2} which are proved by a different method.

We remark that we also have 
$$\kappa(V)\ge\kappa(A, \widehat{\det} f_* \mathcal{O}_X)$$
by using the same method as in the proof of Theorem \ref{m1}, which is the first part of the inequality in Corollary \ref{ineq1} when $m=0$.

Theorem \ref{eq1} and Corollary \ref{ineq1} have many applications, see \cite[Corollaries 1.2, 1.3, 1.7 and 1.8]{Men22}.

With more work, we can generalize Corollary \ref{ineq1} to the case when the base is a smooth projective variety of maximal Albanese dimension. A smooth projective variety is of \emph{maximal Albanese dimension} if its Albanese morphism is generically finite onto its image.

\begin{thm}\label{ineqmad}
Let $f\colon X \to Y$ be a fibration from a smooth projective variety $X$ to a smooth projective variety $Y$ of maximal Albanese dimension where $f$ is smooth over an open set $V\subseteq Y$, and $m$ a positive integer. Then
$$\kappa(V)\ge\kappa(Y, \widehat{\det} f_* \omega_{X}^{\otimes m}\otimes\omega_Y).$$
\end{thm}

By \cite[Theorem 1.1(i)]{Kaw85a}, Theorem \ref{ineqmad} implies a special case of the Kebekus--Kov\'acs conjecture \cite[Conjecture 1.6]{KK08} (see also \cite{PS17}) when the base $V$ compactifies to a smooth projective variety $Y$ of maximal Albanese dimension. This conjecture bounds $\Var(f)$ from above by the log Kodaira dimension $\kappa (V)$ assuming that the general fiber of $f$ has a good minimal model. Progress has been made towards this conjecture, see e.g. \cite{KK08}, \cite{Taj16}, \cite{PS17}.

Next, we consider another interesting problem. Let $f\colon X \to A$ be a surjective morphism from a klt pair $(X, \Delta)$ to an abelian variety $A$, $m\ge1$ a rational number, and $D$ a Cartier divisor on $X$ such that $D\sim_{\Q}m(K_X+\Delta)$. We decide when
$$\kappa(A, \widehat{\det} f_*\mathcal{O}_X(D))=\dim V^0(A, f_*\mathcal{O}_X(D))$$
holds using the techniques of this paper. These quantities have appeared many times in this paper and are very interesting in their own right. We have that $\kappa(A, \widehat{\det} f_*\mathcal{O}_X(D))\ge\dim V^0(A, f_*\mathcal{O}_X(D))$, and $\kappa(A, \widehat{\det} f_*\mathcal{O}_X(D))=\dim V^0(A, f_*\mathcal{O}_X(D))$ if $m>1$ by \cite[Theorem 3.6, Remark 3.7 and Lemma 3.8]{Men22}. However, the inequality can be strict if $m=1$. For example, in \cite[Example 1.13]{EL97}, Ein and Lazarsfeld constructed a smooth projective threefold $X$ of general type and of maximal Albanese dimension and a generically finite surjective morphism $f\colon X \to A$ to an abelian threefold $A$ such that $\kappa(A, \widehat{\det} f_* \omega_{X})=3>\dim V^0(A, f_*\omega_X)=2$, see also \cite[Example 4.1]{CDJ14}. We decide when the equality holds as an application of Theorem \ref{eqcon}.

\begin{thm}\label{eqconc}
Let $f\colon X \to A$ be a surjective morphism from a klt pair $(X, \Delta)$ to an abelian variety $A$, $m\ge1$ a rational number, and $D$ a Cartier divisor on $X$ such that $D\sim_{\Q}m(K_X+\Delta)$. Then
$$\kappa(A, \widehat{\det} f_*\mathcal{O}_X(D))=\dim V^0(A, f_*\mathcal{O}_X(D))$$
if and only if $V^0(A, f_*\mathcal{O}_X(D))^{\otimes M}$ is irreducible for sufficiently divisible positive integer $M$.
\end{thm}

The set $V^0(A, f_*\mathcal{O}_X(D))^{\otimes M}$ is defined by
$$V^0(A, f_*\mathcal{O}_X(D))^{\otimes M}=\{\alpha^{\otimes M}\in\Pic^0(A)\mid\alpha\in V^0(A, f_*\mathcal{O}_X(D))\},$$
which is a closed subset of $\Pic^0(A)$.

If $m>1$, it is known that $V^0(A, f_*\mathcal{O}_X(D))^{\otimes M}$ is irreducible for sufficiently divisible positive integer $M$ by \cite[Theorem 2.6]{Men22} and \cite[Lemma 3.3]{LPS20}, and $\kappa(A, \widehat{\det} f_*\mathcal{O}_X(D))=\dim V^0(A, f_*\mathcal{O}_X(D))$ by \cite[Lemma 3.8]{Men22}. Thus Theorem \ref{eqconc} follows automatically if $m>1$. However, things become subtle if $m=1$. In the previous example constructed by Ein and Lazarsfeld, $V^0(A, f_*\omega_X)^{\otimes M}$ has three irreducible components for every positive integer $M$, and $\kappa(A, \widehat{\det} f_*\omega_{X})>\dim V^0(A, f_*\omega_X)$.

\subsection{A generalization of the Iitaka conjecture}

In this subsection, we present more applications. We consider an Iitaka-type conjecture and give evidence. Iitaka proposed the following conjecture which predicts the subadditivity of Kodaira dimension of fibrations.

\begin{Iitaka conjecture}
Let $f\colon X \to Y$ be a fibration between smooth projective varieties with general fiber $F$. Then
$$\kappa(X)\ge\kappa(F)+\kappa(Y).$$
\end{Iitaka conjecture}

The Iitaka conjecture is known in various cases but not fully known yet. For example, it is known when $Y$ is of general type by \cite{Vie83a}, $Y$ is a curve by \cite{Kaw82}, $F$ has a good minimal model by \cite{Kaw85a}, $F$ is of general type by \cite{Kol87}, $Y$ is an abelian variety by \cite{CP17} (see also \cite{HPS18}), $Y$ is of maximal Albanese dimension by \cite{HPS18}, or $Y$ is a surface by \cite{Cao18}. The logarithmic version of this conjecture is also known in several cases, see e.g. \cite{Fuj17} and \cite{KP17}.

Viehweg proposed the following conjecture as a strengthening of the Iitaka conjecture.

\begin{Viehweg conjecture}
Let $f\colon X \to Y$ be a fibration between smooth projective varieties with general fiber $F$ and $\kappa(Y)\ge0$. Then,
$$\kappa(X)\ge\kappa(F)+\max\{\kappa(Y),\,\Var(f)\}.$$
\end{Viehweg conjecture}

Unlike the Iitaka conjecture, the Viehweg conjecture is only known in very few cases. For example, it is known when $Y$ is of general type by \cite{Vie83a}, $F$ has a good minimal model by \cite{Kaw85a}, or $F$ is of general type by \cite{Kol87}.

Recently, Popa proposed several conjectures on the superadditivity of Kodaira dimension of fibrations, see \cite{Pop21} for details.

In a different but related direction, we propose the following conjecture which predicts an equality of Kodaira dimension of fibrations.

\begin{conj}\label{sic}
Let $f\colon X \to Y$ be a fibration between smooth projective varieties with general fiber $F$ and $\kappa(Y)\ge0$. For sufficiently divisible positive integer $m$ such that $f_* \omega_{X}^{\otimes m}\neq0$, then
$$\kappa(X)=\kappa(F)+\kappa(Y, \widehat{\det} f_* \omega_{X}^{\otimes m}\otimes\omega_Y).$$
\end{conj}

If $\kappa(F)=-\infty$, then $\kappa(X)=-\infty$ and the equality is trivial. If $\kappa(F)\ge0$, the good minimal model conjecture implies that $F$ has a good minimal model. We deduce that
$$\kappa(Y, \widehat{\det} f_* \omega_{X/Y}^{\otimes m})\ge\Var(f)\ge0$$
for some sufficiently divisible positive integer $m$ by \cite[Theorem 1.1(i)]{Kaw85a}. Since $\kappa(Y)\ge0$, we deduce that
$$\kappa(Y, \widehat{\det} f_* \omega_{X}^{\otimes m}\otimes\omega_Y)=\kappa(Y, \widehat{\det} f_* \omega_{X/Y}^{\otimes m}\otimes\omega_Y^{\otimes(mr+1)})\ge\max\{\kappa(Y),\,\Var(f)\}$$
where $r=\rank f_* \omega_{X/Y}^{\otimes m}$. Thus Conjecture \ref{sic} implies Viehweg conjecture assuming the good minimal model conjecture.

We prove that Conjecture \ref{sic} is true in several cases.

\begin{thm}\label{sicpr}
Conjecture \ref{sic} is true when
\begin{enumerate}
	\item[(1)] $Y$ is of maximal Albanese dimension.

	\item[(2)] $Y$ is of general type.
	
	\item[(3)] $Y$ is a curve.
\end{enumerate}
\end{thm}

The case when $Y$ is of general type follows directly from \cite[Theorem III and Corollary IV]{Vie83a}. The case when $Y$ is a curve follows from the case when $Y$ is of maximal Albanese dimension since a curve with nonnegative Kodaira dimension is of maximal Albanese dimension. For the case when $Y$ is of maximal Albanese dimension, we prove the following stronger theorem. Theorem \ref{eq1} is an important ingredient in its proof.

\begin{thm}\label{eqmad}
Let $f\colon X \to Y$ be a fibration from a klt pair $(X, \Delta)$ to a smooth projective variety $Y$ of maximal Albanese dimension, $F$ the general fiber of $f$, $m>1$ a rational number, and $D$ a Cartier divisor on $X$ such that $D\sim_{\Q}m(K_X+\Delta)$ and $f_*\mathcal{O}_X(D)\neq0$. Let $L$ be a divisor on $Y$ such that $\kappa(Y, L)\ge0$. For every positive integer $l$, we have
$$\kappa(X, K_X+\Delta+f^*L)=\kappa(F, K_F+\Delta|_F)+\kappa(Y, \widehat{\det} f_*\mathcal{O}_X(lD)\otimes\mathcal{O}_Y(K_Y+L)).$$
\end{thm}

\subsection{Main theorems}

Let $f\colon X \to A$ be a surjective morphism from a smooth projective variety $X$ to an abelian variety $A$. Our first main theorem shows that there is a connection between the singularity of $f$ and the positivity of $\widehat{\det} f_* \omega_{X}^{\otimes m}$.

\begin{thm}\label{m1}
Let $f\colon X \to A$ be a surjective morphism from a smooth projective variety $X$ to an abelian variety $A$ where $f$ is smooth away from a closed set $Z\subseteq A$ such that $Z=W\cup D$ where $W\subseteq A$ is a closed set with $\codim_A W\ge 2$ and $D$ is an effective divisor on $A$. Let $m$ be a positive integer and $T$ a divisor on $A$ such that $f_* \omega_{X}^{\otimes m}\neq0$ and $\mathcal{O}_A(T)\cong\widehat{\det} f_* \omega_{X}^{\otimes m}$. Then
\begin{enumerate}

	\item[(1)] There exists a fibration $h\colon A\to C$ to an abelian variety $C$ such that $T\sim_{\Q}h^*N$ where $N$ is an ample effective $\Q$-divisor on $C$.

	\item[(2)] Let $h\colon A\to C$ be any fibration satisfying statement (1). Let $g\colon A\to B$ be any fibration to an abelian variety B such that $D=g^*H$ where $H$ is an ample effective divisor on $B$. Then there exists a fibration $t\colon B\to C$ such that $h=t\circ g$.
\end{enumerate}
\[
	\begin{tikzcd}
			X \dar{f} \\
			A \dar{g} \rar{h} & C \\
			B \urar[swap]{t}
	\end{tikzcd}
\]
\end{thm}

The fibration $g\colon A\to B$ in the theorem above always exists due to a well-known structural theorem for effective divisors on abelian varieties.

Statement (1) of the theorem above follows directly from \cite[Lemma 3.3]{Men22} and a well-known structural theorem for effective divisors on abelian varieties. The central part of Theorem \ref{m1} is its statement (2). The positivity of the line bundle $\widehat{\det} f_* \omega_{X}^{\otimes m}$ is encoded by the abelian variety $C$ and the ample effective $\Q$-divisor $N$ on $C$. Similarly, the positivity of the divisor $D$ is encoded by the abelian variety $B$ and the ample effective divisor $H$ on $B$. Thus statement (2) of the theorem above shows that the positivity of $\widehat{\det} f_* \omega_{X}^{\otimes m}$ is controlled by the positivity of $D$ in a concrete way which is given by the fibration $t\colon B\to C$.

The next theorem compares the positivity of the line bundle $\widehat{\det} f_* \omega_{X}^{\otimes m}$ with the positivity of the coherent sheaf $f_* \omega_{X}^{\otimes m}$. The positivity of $f_* \omega_{X}^{\otimes m}$ is described by its Chen--Jiang decomposition, see Definition \ref{CJD}. It is known that pushforwards of pluricanonical bundles under morphisms to abelian varieties have the Chen--Jiang decomposition by \cite{CJ18, PPS17, LPS20} in increasing generality, and pushforwards of klt pairs under morphisms to abelian varieties have the Chen--Jiang decomposition, as proved independently in \cite{Jia22} and \cite{Men21}. 

\begin{thm}\label{m2}
Let $f\colon X \to A$ be a surjective morphism from a klt pair $(X, \Delta)$ to an abelian variety $A$, $m\geq1$ a rational number, and $D$ a Cartier divisor on $X$ such that $D\sim_{\Q}m(K_X+\Delta)$ and $f_*\mathcal{O}_X(D)\neq0$. Let $T$ be a divisor on $A$ such that $\mathcal{O}_A(T)\cong\widehat{\det} f_*\mathcal{O}_X(D)$. Then
\begin{enumerate}

	\item[(1)] There exists a fibration $h\colon A\to C$ to an abelian variety $C$ such that $T\sim_{\Q}h^*N$ where $N$ is an ample effective $\Q$-divisor on $C$.

	\item[(2)] Let $h\colon A\to C$ be any fibration satisfying statement (1). Let
	$$f_*\mathcal{O}_X(D)\cong \bigoplus_{i\in I}(\alpha_i\otimes p_i^*\mathcal{F}_i)$$
	be any finite direct sum decomposition of $f_*\mathcal{O}_X(D)$ such that each $A_i$ is an abelian variety, each $p_i\colon A\to A_i$ is a fibration, each $\mathcal{F}_i$ is a nonzero M-regular coherent sheaf on $A_i$, and each $\alpha_i\in\Pic^0(A)$ is a torsion line bundle. Then there exists a fibration $t_i\colon C\to A_i$ such that $p_i=t_i\circ h$ for every $i\in I$.
\end{enumerate}
\[
	\begin{tikzcd}
			(X, \Delta) \dar{f} \\
			A \dar{h} \rar{p_i} & A_i \\
			C \urar[swap]{t_i}
	\end{tikzcd}
\]
\end{thm}

Under the hypotheses of Theorem \ref{m2}, the decomposition of $f_*\mathcal{O}_X(D)$ as in statement (2) always exists by \cite[Theorem 1.3]{Men21}, see also \cite[Theorem 1.3]{Jia22} for the case when $m\geq1$ is an integer. The theorem shows that each fibration $p_i$ arising from any Chen--Jiang decomposition of $f_*\mathcal{O}_X(D)$ factors through $h$.

Statement (1) of the theorem above follows directly from \cite[Lemma 3.3]{Men22} and a well-known structural theorem for effective divisors on abelian varieties. The central part of Theorem \ref{m2} is its statement (2), which shows that the positivity of $f_*\mathcal{O}_X(D)$ is controlled by the positivity of $\widehat{\det} f_*\mathcal{O}_X(D)$ in a concrete way which is given by the fibrations $t_i\colon C\to A_i$. The proof of Theorem \ref{m2} is more complicated than the proof of Theorem \ref{m1} and requires more techniques.

We immediately have the following corollary of Theorems \ref{m1} and \ref{m2} by combining them together.

\begin{coro}\label{rel}
Let $f\colon X \to A$ be a surjective morphism from a smooth projective variety $X$ to an abelian variety $A$ where $f$ is smooth away from a closed set $Z\subseteq A$ such that $Z=W\cup D$ where $W\subseteq A$ is a closed set with $\codim_A W\ge 2$ and $D$ is an effective divisor on $A$. Let $g\colon A\to B$ be any fibration to an abelian variety B such that $D=g^*H$ where $H$ is an ample effective divisor on $B$. Let $m$ be any positive integer such that $f_* \omega_{X}^{\otimes m}\neq0$. Let
$$f_* \omega_{X}^{\otimes m}\cong \bigoplus_{i\in I}(\alpha_i\otimes p_i^*\mathcal{F}_i)$$
be any finite direct sum decomposition of $f_* \omega_{X}^{\otimes m}$ such that each $A_i$ is an abelian variety, each $p_i\colon A\to A_i$ is a fibration, each $\mathcal{F}_i$ is a nonzero M-regular coherent sheaf on $A_i$, and each $\alpha_i\in\Pic^0(A)$ is a torsion line bundle. Then there exists a fibration $s_i\colon B\to A_i$ such that $p_i=s_i\circ g$ for every $i\in I$. 

\[
	\begin{tikzcd}
			X \dar{f} \\
			A \dar{g} \rar{p_i} & A_i \\
			B \urar[swap]{s_i}
	\end{tikzcd}
\]
\end{coro}

Corollary \ref{rel} establishes relationship between the singularity of the surjective morphism $f$ and the positivity of $f_* \omega_{X}^{\otimes m}$. Concretely speaking, it establishes relationship between the divisorial part $D$ of the closed set $Z\subseteq A$ which $f$ is smooth away from and the Chen--Jiang decomposition of $f_* \omega_{X}^{\otimes m}$. It can be viewed as a generalization of \cite[Theorem B]{MP23} to the case when $D$ is not necessarily the zero divisor. \cite[Theorem B]{MP23} gives the structures of pushforwards of pluricanonical bundles of smooth projective varieties under surjective morphisms to abelian varieties when the morphisms are smooth away from a closed set of codimension at least $2$ in the abelian varieties.

In the proofs of the main theorems, we employ results from \cite{Men21}, \cite{MP23} and \cite{Men22}, and arguments on positivity properties of coherent sheaves. The line bundle $\widehat{\det} f_*\mathcal{O}_X(D)$ plays an important role in the proofs.

\begin{ac}
{I would like to thank Jungkai Alfred Chen, Christopher D. Hacon, Mircea Musta{\c{t}}{\u{a}}, Mihnea Popa, Vyacheslav Shokurov and Ziquan Zhuang for helpful discussions. I would also like to thank the referees for helpful comments.}
\end{ac}

\section{Preliminaries}\label{2}

We work over $\C$. Let $\mathcal{F}$ be a coherent sheaf on a projective variety $X$. We denote $\mathcal{H}om_{\mathcal{O}_{X}}(\mathcal{F}, \mathcal{O}_{X})$ by $\mathcal{F}^{\vee}$. Let $L$ be a $\Q$-Cartier $\Q$-divisor on a normal projective variety $X$. We denote the \emph{Iitaka dimension} of $L$ by $\kappa (X, L)$, see \cite{Kaw85b}. A \emph{fibration} is a projective surjective morphism with connected fibers.

We recall several definitions first.

\begin{defi}\label{irregularity}
Let $X$ be a smooth projective variety. The \emph{irregularity} $q(X)$ is defined as $h^1(X,\mathcal{O}_X)$. If $X$ is a projective variety, the \emph{irregularity} $q(X)$ is defined as the irregularity of any resolution of $X$.
\end{defi}

Note that the irregularity of a projective variety is well-defined. If $X$ is a normal projective variety with rational singularities, then the irregularity $q(X)$ is equal to the dimension of its Albanese variety $\Alb(X)$ since its Albanese variety coincides with the Albanese variety of any of its resolutions by \cite[Proposition 2.3]{Rei83} and \cite[Lemma 8.1]{Kaw85a}.

\begin{defi}\label{csl}
Let $\mathcal{F}$ be a coherent sheaf on an abelian variety $A$. The \emph{cohomological support loci} $V_l^i(A, \mathcal{F})$ for $i\in\N$ and $l\in\N$ are defined by
$$V_l^i(A, \mathcal{F})=\{\alpha\in\Pic^0(A)\mid\dim H^i(A, \mathcal{F}\otimes\alpha)\geq l\}.$$
We use $V^i(A, \mathcal{F})$ to denote $V_1^i(A, \mathcal{F})$.
\end{defi}

\begin{defi}\label{GVM}
A coherent sheaf $\mathcal{F}$ on an abelian variety $A$
\begin{enumerate}
	\item[(1)] is a GV-\emph{sheaf} if $\codim_{\Pic^0 (A)} V^i (A, \mathcal{F}) \ge i$ for every $i>0$.
	\item[(2)] is \emph{M-regular} if $\codim_{\Pic^0 (A)} V^i (A, \mathcal{F}) > i$ for every $i>0$.
	\item[(3)] \emph{satisfies} $\IT$ if $V^i (A, \mathcal{F})=\emptyset$ for every $i>0$.
\end{enumerate}
\end{defi}

It is known that M-regular sheaves are ample by \cite[Corollary 3.2]{Deb06}, and GV-sheaves are nef by \cite[Theorem 4.1]{PP11b}. 

We now give the definition of the Chen--Jiang decomposition.

\begin{defi}\label{CJD}
Let $\mathcal{F}$ be a coherent sheaf on an abelian variety $A$. The sheaf $\mathcal{F}$ is said to have the \emph{Chen--Jiang decomposition} if $\mathcal{F}$ admits a finite direct sum decomposition
$$\mathcal{F}\cong \bigoplus_{i\in I}(\alpha_i\otimes p_i^*\mathcal{F}_i),$$
where each $A_i$ is an abelian variety, each $p_i\colon A\to A_i$ is a fibration, each $\mathcal{F}_i$ is a nonzero M-regular coherent sheaf on $A_i$, and each $\alpha_i\in\Pic^0(A)$ is a torsion line bundle.
\end{defi}

\section{Proofs of main results}\label{3}

We prove the main theorems first. We start with Theorem \ref{m1}.

\begin{proof}[Proof of Theorem \ref{m1}]
We prove statement (1) first. By \cite[Lemma 3.3]{Men22}, we have
$$\kappa(A, \mathcal{O}_A(T))=\kappa(A, \widehat{\det} f_* \omega_{X}^{\otimes m})\ge0.$$
Thus there exists a positive integer $k$ such that $kT\sim E\ge0$ where $E$ is an effective divisor on $A$. By a well-known structural theorem, there exist a fibration $h\colon A \to C$ between abelian varieties and an ample effective divisor $N'$ on $C$ such that $E=h^*N'$. We let $N:=\frac{N'}{k}$. Then $N$ is an ample effective $\Q$-divisor on $C$, and we have
$$T\sim_{\Q}\frac{E}{k}=\frac{h^*N'}{k}=h^*N.$$

Next, we prove statement (2). Let $h\colon A\to C$ be a fibration satisfying statement (1). Let $g\colon A\to B$ be a fibration to an abelian variety $B$ such that $D=g^*H$ where $H$ is an ample effective divisor on $B$. We consider the following commutative diagram where $b\in B$ is a general closed point, and $f_{b}$ is obtained by base change from $f$ via the closed immersion $i_b$ of the closed fiber $A_b$ of $g$ over $b$.
\[
	\begin{tikzcd}
	X_b\dar{f_b} \rar	&	X \dar{f} \\
	A_b\dar \rar{i_b}	&	A \dar{g} \rar{h} & C \\
	\{b\} \rar	&	B
	\end{tikzcd}
\]
We can choose $b$ sufficiently general such that $X_b$ is smooth, $f_b$ is surjective, $i_b^{-1}(D)=\emptyset$, and $\codim_{A_b}  i_b^{-1}(W)\ge 2$ since $D=g^*H$ and $\codim_A W\ge 2$. By \cite[Lemma 3.2]{Men22}, we can choose $b$ sufficiently general such that
$$i_b^*\mathcal{O}_A(T)\cong i_b^*\widehat{\det} f_* \omega_{X}^{\otimes m}\cong\widehat{\det} {f_b}_*(\omega_{X}^{\otimes m}|_{X_b})\cong \widehat{\det} {f_b}_*\omega_{X_b}^{\otimes m}.$$
Since $\codim_{A_b}  i_b^{-1}(W)\ge 2$ and $f_b$ is surjective and smooth away from $i_b^{-1}(Z)=i_b^{-1}(W)$, we deduce that
$${f_b}_*\omega_{X_b}^{\otimes m}\cong\bigoplus_{i=1}^{s}\beta_i$$
by \cite[Theorem B]{MP23} where $s$ is a positive integer and each $\beta_i$ is a torsion line bundle. Thus
$$i_b^*\mathcal{O}_A(T)\cong\widehat{\det} {f_b}_*\omega_{X_b}^{\otimes m}\cong\bigotimes_{i=1}^{s}\beta_i $$
is a torsion line bundle. In particular, we have $\kappa(A_b, i_b^*T)=0$.

We claim that the image of $A_b$ under the morphism $h\circ i_b$ is a point for sufficiently general closed point $b$. We prove the claim by contradiction. Assume that it is not a point. Then the scheme theoretic image denoted by $(h\circ i_b)(A_b)$ of $h\circ i_b$ is an abelian variety of positive dimension since $h\circ i_b$ is a morphism between abelian varieties. We consider the following commutative diagram where $v$ is the closed immersion from $(h\circ i_b)(A_b)$ to $C$, and $u$ is the surjective morphism induced by $h\circ i_b$.
\[
	\begin{tikzcd}
	A_b\dar{u} \rar{i_b}	&	A \dar{h} \\
	(h\circ i_b)(A_b)\rar{v} &	C
	\end{tikzcd}
\]
Since $v$ is a closed immersion, we deduce that $v^*N$ is an ample $\Q$-divisor on $(h\circ i_b)(A_b)$. Since $u$ is surjective, we have
$$0=\kappa(A_b, i_b^*T)=\kappa(A_b, i_b^*h^*N)=\kappa(A_b, u^*v^*N)$$
$$=\kappa((h\circ i_b)(A_b), v^*N)=\dim (h\circ i_b)(A_b)>0,$$
which is a contradiction. Thus $(h\circ i_b)(A_b)$ must be a point, and we finish proving the claim. Since $B$ is smooth and $C$ is an abelian variety, we deduce that there exists a morphism $t\colon B\to C$ such that $h=t\circ g$ by \cite[Lemma 14]{Kaw81}. Since $g$ and $h$ are fibrations, $t$ is also a fibration.
\end{proof}

Next, we prove Theorem \ref{m2}, which requires more techniques.

\begin{proof}[Proof of Theorem \ref{m2}]
Statement (1) follows from the same argument as in the proof of Theorem \ref{m1}.

We prove statement (2) now. Let $h\colon A\to C$ be a fibration satisfying statement (1). Let
	$$f_*\mathcal{O}_X(D)\cong \bigoplus_{i\in I}(\alpha_i\otimes p_i^*\mathcal{F}_i)$$
	be a finite direct sum decomposition of $f_*\mathcal{O}_X(D)$ such that each $A_i$ is an abelian variety, each $p_i\colon A\to A_i$ is a fibration, each $\mathcal{F}_i$ is a nonzero M-regular coherent sheaf on $A_i$, and each $\alpha_i\in\Pic^0(A)$ is a torsion line bundle. We consider the following commutative diagram where $c\in C$ is a general closed point, and $f_{c}$ is obtained by base change from $f$ via the closed immersion $i_c$ of the closed fiber $A_c$ of $h$ over $c$.
\[
	\begin{tikzcd}
	(X_c, \Delta|_{X_c})\dar{f_c} \rar	&	(X, \Delta) \dar{f} \\
	A_c\dar \rar{i_c}	&	A \dar{h} \rar{p_i} & A_i \\
	\{c\} \rar	&	C
	\end{tikzcd}
\]
We can choose $c$ sufficiently general such that $(X_c, \Delta|_{X_c})$ is klt, and $f_c$ is surjective. We have
$$\kappa(A_c, i_c^*\widehat{\det} f_*\mathcal{O}_X(D))=\kappa(A_c, i_c^*T)=\kappa(A_c, i_c^*h^*N)=0.$$
We have that $\mathcal{F}_i$ is a direct summand of ${p_i}_*f_*(\mathcal{O}_X(D)\otimes f^*\alpha_i^{-1})$ for every $i\in I$. Thus we deduce that $\mathcal{F}_i$ is torsion-free since $p_i\circ f$ is surjective. Since $\mathcal{F}_i$ is an M-regular sheaf on $A_i$, it is ample by \cite[Proposition 2.13]{PP03} and \cite[Corollary 3.2]{Deb06}. The sheaf $\mathcal{F}_i$ is big since an ample sheaf is big (see e.g. \cite[Section 2]{Deb06} and \cite[Section 5]{Mor87}). Thus $\widehat{\det}\mathcal{F}_i$ is a big line bundle by \cite[Lemma 3.2(iii)]{Vie83b} (see also \cite[Properties 5.1.1]{Mor87}). Since $A_i$ is an abelian variety, $\widehat{\det}\mathcal{F}_i$ is an ample line bundle. Since $p_i$ is flat and $\widehat{\det}(\alpha_i\otimes p_i^*\mathcal{F}_i)$ is a line bundle, we deduce that
$$\widehat{\det} f_*\mathcal{O}_X(D)\cong\bigotimes_{i\in I}\widehat{\det}(\alpha_i\otimes p_i^*\mathcal{F}_i)\cong\bigotimes_{i\in I}(\widehat{\det}p_i^*\mathcal{F}_i\otimes\alpha_i^{\otimes\rank\mathcal{F}_i})$$
$$\cong\bigotimes_{i\in I}(p_i^*\widehat{\det}\mathcal{F}_i\otimes\alpha_i^{\otimes\rank\mathcal{F}_i})\cong (\bigotimes_{i\in I}p_i^*\widehat{\det}\mathcal{F}_i)\otimes(\bigotimes_{i\in I}\alpha_i^{\otimes\rank\mathcal{F}_i}).$$
Since $\widehat{\det}\mathcal{F}_i$ is an ample line bundle, we deduce that $\kappa(A_c, i_c^*p_i^*\widehat{\det}\mathcal{F}_i)\ge0$. Thus we deduce that
$$0=\kappa(A_c, i_c^*\widehat{\det} f_*\mathcal{O}_X(D))=\kappa(A_c, (\bigotimes_{i\in I}i_c^*p_i^*\widehat{\det}\mathcal{F}_i)\otimes(\bigotimes_{i\in I}i_c^*\alpha_i^{\otimes\rank\mathcal{F}_i}))$$
$$=\kappa(A_c, \bigotimes_{i\in I}i_c^*p_i^*\widehat{\det}\mathcal{F}_i)\ge\kappa(A_c, i_c^*p_i^*\widehat{\det}\mathcal{F}_i)\ge0.$$
Thus for every $i\in I$, we deduce that 
$$\kappa(A_c, i_c^*p_i^*\widehat{\det}\mathcal{F}_i)=0.$$

We claim that the image of $A_c$ under the morphism $p_i\circ i_c$ is a point for sufficiently general closed point $c$. We prove the claim by contradiction. Assume that it is not a point. Then the scheme theoretic image denoted by $(p_i\circ i_c)(A_c)$ of $p_i\circ i_c$ is an abelian variety of positive dimension since $p_i\circ i_c$ is a morphism between abelian varieties. We consider the following commutative diagram where $v$ is the closed immersion from $(p_i\circ i_c)(A_c)$ to $A_i$, and $u$ is the surjective morphism induced by $p_i\circ i_c$.
\[
	\begin{tikzcd}
	A_c\dar{u} \rar{i_c}	&	A \dar{p_i} \\
	(p_i\circ i_c)(A_c)\rar{v} &	A_i
	\end{tikzcd}
\]
Since $v$ is a closed immersion, we deduce that $v^*\widehat{\det}\mathcal{F}_i$ is an ample line bundle on $(p_i\circ i_c)(A_c)$. Since $u$ is surjective, we have
$$0=\kappa(A_c, i_c^*p_i^*\widehat{\det}\mathcal{F}_i)=\kappa(A_c, u^*v^*\widehat{\det}\mathcal{F}_i)$$
$$=\kappa((p_i\circ i_c)(A_c), v^*\widehat{\det}\mathcal{F}_i)=\dim (p_i\circ i_c)(A_c)>0,$$
which is a contradiction. Thus $(p_i\circ i_c)(A_c)$ must be a point, and we finish proving the claim. Since $C$ is smooth and $A_i$ is an abelian variety, we deduce that there exists a morphism $t_i\colon C\to A_i$ such that $p_i=t_i\circ h$ for every $i\in I$ by \cite[Lemma 14]{Kaw81}. Since $h$ and $p_i$ are fibrations, $t_i$ is also a fibration.
\end{proof}


\begin{proof}[Proof of Corollary \ref{ineq1}]
If $f_* \omega_{X}^{\otimes m}=0$, then the statement is trivial. Thus we can assume that $f_* \omega_{X}^{\otimes m}\neq0$. By \cite[Theorem 3.6, Remark 3.7 and Lemma 3.8]{Men22}, we only need to prove that
$$\kappa(V)\ge\kappa(A, \widehat{\det} f_* \omega_{X}^{\otimes m}).$$
Let $T$ be a divisor on $A$ such that $\mathcal{O}_A(T)\cong\widehat{\det} f_* \omega_{X}^{\otimes m}$. Let $Z:=A\setminus V$ which is a closed subset of $A$. Then $Z=W\cup D$ where $W\subseteq A$ is a closed set with $\codim_A W\ge 2$ and $D$ is an effective divisor on $A$. By a well-known structural theorem, there exist a fibration $g\colon A \to B$ between abelian varieties and an ample effective divisor $H$ on $B$ such that $D=g^*H$. By \cite[Lemma 2.6]{MP23}, we have that
$$\kappa(V)=\kappa(A, K_A+D)=\kappa(A, g^*H)=\kappa(B, H)=\dim B.$$
By Theorem \ref{m1}, there exist a fibration $h\colon A\to C$ to an abelian variety $C$ such that $T\sim_{\Q}h^*N$ where $N$ is an ample effective $\Q$-divisor on $C$ and a fibration $t\colon B\to C$ such that $h=t\circ g$. Thus
$$\kappa(V)=\dim B\ge\dim C=\kappa(C, N)$$
$$=\kappa(A, h^*N)=\kappa(A, T)=\kappa(A, \widehat{\det} f_* \omega_{X}^{\otimes m}).$$
\end{proof}

Next, we prove Theorem \ref{eq1} which gives an equality of Kodaira dimension of surjective morphisms to abelian varieties.

\begin{proof}[Proof of Theorem \ref{eq1}]
We only need to prove the equality when $l=1$. By \cite[Theorem 2.6]{Men22}, there exists a fibration $p\colon A\to B$ to an abelian variety $B$ such that $f_*\mathcal{O}_X(kD)$ admits, for every positive integer $k$, a finite direct sum decomposition
$$f_*\mathcal{O}_X(kD)\cong\bigoplus_{i\in I}(\alpha_i\otimes p^*\mathcal{F}_i),$$
where each $\mathcal{F}_i$ is a nonzero coherent sheaf on $B$ satisfying $\IT$, and each $\alpha_i\in\Pic^0(A)$ is a torsion line bundle. We consider the following commutative diagram where $b\in B$ is a general closed point, and $f_{b}$ is obtained by base change from $f$ via the closed immersion $i_b$ of the closed fiber $A_b$ of $p$ over $b$.
\[
	\begin{tikzcd}
	(X_b, \Delta|_{X_b})\dar{f_b} \rar	&	(X, \Delta) \dar{f} \\
	A_b\dar \rar{i_b}	&	A \dar{p} \\
	\{b\} \rar	&	B
	\end{tikzcd}
\]
We can choose $b$ sufficiently general such that $(X_b, \Delta|_{X_b})$ is klt and $f_b$ is a surjective morphism. We have that $\mathcal{F}_i$ is a direct summand of $p_*f_*(\mathcal{O}_X(kD)\otimes f^*\alpha_i^{-1})$. Since $\mathcal{F}_i$ satisfies $\IT$, $\mathcal{F}_i$ is M-regular and thus ample by \cite[Proposition 2.13]{PP03} and \cite[Corollary 3.2]{Deb06}. By \cite[Lemma 3.4]{Men22} which works for surjective morphisms by the same argument as in the proof of itself, we deduce that
$$\kappa(X, D)=\kappa(X, \mathcal{O}_X(kD)\otimes f^*\alpha_i^{-1})$$
$$=\kappa(X_b, \mathcal{O}_{X_b}(kD|_{X_b})\otimes (f^*\alpha_i^{-1})|_{X_b})+\dim B=\kappa(X_b, D|_{X_b})+\dim B$$
for very general $b$. We can choose very general $b$ such that
$${f_b}_*\mathcal{O}_{X_b}(kD|_{X_b})\cong i_b^*f_*\mathcal{O}_{X}(kD)$$
for every positive integer $k$ by \cite[Proposition 4.1]{LPS20}. Thus we deduce that
$${f_b}_*\mathcal{O}_{X_b}(kD|_{X_b})\cong i_b^*f_*\mathcal{O}_{X}(kD)\cong\bigoplus_{i\in I}(i_b^*\alpha_i\otimes i_b^*p^*\mathcal{F}_i).$$
Thus ${f_b}_*\mathcal{O}_{X_b}(kD|_{X_b})$ is a direct sum of torsion line bundles for every positive integer $k$. The general fiber of $f_b$ is $F$. We deduce that
$$h^0(X_b, \mathcal{O}_{X_b}(kD|_{X_b}))=h^0(A_b, {f_b}_*\mathcal{O}_{X_b}(kD|_{X_b}))$$
$$\le\rank {f_b}_*\mathcal{O}_{X_b}(kD|_{X_b})\le h^0(F, \mathcal{O}_{F}(kD|_{F}))$$
for every positive integer $k$. Thus $\kappa(X_b, D|_{X_b})\le\kappa(F, D|_{F})$. By \cite[Theorem 1.1]{Men21}, we have $\kappa(X_b, D|_{X_b})\ge\kappa(F, D|_{F})$ since $A_b$ is an abelian variety. Thus
$$\kappa(X_b, D|_{X_b})=\kappa(F, D|_F).$$
By \cite[Theorem 4.2]{HMX18}, $\kappa(F, D|_F)$ is constant for general fiber $F$ of $f$ since $(X, \Delta)$ is klt. We deduce that 
$$\kappa(X, K_X+\Delta)=\kappa(X, D)=\kappa(X_b, D|_{X_b})+\dim B$$
$$=\kappa(F, D|_F)+\dim B=\kappa(F, K_F+\Delta|_F)+\dim B.$$
By \cite[Lemma 3.3]{LPS20}, we deduce that
$$\dim V^0(A, f_*\mathcal{O}_X(D))=\dim B.$$
By \cite[Lemma 3.8]{Men22}, we have that
$$\kappa(A, \widehat{\det} f_*\mathcal{O}_X(D))=\dim V^0(A, f_*\mathcal{O}_X(D)).$$
Thus we have
$$\kappa(X, K_X+\Delta)=\kappa(F, K_F+\Delta|_F)+\dim V^0(A, f_*\mathcal{O}_X(D))$$
$$=\kappa(F, K_F+\Delta|_F)+\kappa(A, \widehat{\det} f_*\mathcal{O}_X(D)).$$
\end{proof}


\begin{proof}[Proof of Corollary \ref{eq1c2}]
We can assume that $\kappa(F, K_F+\Delta|_F)\ge0$. Thus we can choose a sufficiently divisible positive integer $N>1$ such that $N(K_X+\Delta)$ is a Cartier divisor and $f_*\mathcal{O}_X(N(K_X+\Delta))\neq0$. By \cite[Theorem 2.6]{Men22} and \cite[Lemma 3.3]{LPS20}, there exists an abelian variety $B$ such that
$$\dim V^0(A, f_*\mathcal{O}_X(N(K_X+\Delta)))=\dim B\le q(Y).$$
By Theorem \ref{eq1}, we deduce that
$$\kappa(X, K_X+\Delta)=\kappa(F, K_F+\Delta|_F)+\dim V^0(A, f_*\mathcal{O}_X(N(K_X+\Delta)))$$
$$\le\kappa(F, K_F+\Delta|_F)+q(Y).$$
\end{proof}

\begin{proof}[Proof of Theorem \ref{eq2}]
Since $\kappa(A, L)\ge0$ and $A$ is an abelian variety, $L$ is semiample and thus $f^*L$ is semiample. Thus we can find an effective $\Q$-divisor $E$ such that $E\sim_{\Q}\frac{f^*L}{lm}$ and $(X, \Delta+E)$ is klt. We have that
$$lD+f^*L\sim_{\Q}lm(K_X+\Delta+\frac{f^*L}{lm})\sim_{\Q}lm(K_X+\Delta+E)$$
where $lm>1$ is a rational number. Since $f_*\mathcal{O}_X(D)\neq0$, $\kappa(F, D|_F)\ge0$ and thus $\kappa(X, D)\ge0$ by Theorem \ref{eq1}. Since $\kappa(X, D)\ge0$ and $\kappa(X, f^*L)\ge0$, we deduce that $\kappa(X, K_X+\Delta+f^*L)=\kappa(X, lD+f^*L)$. By Theorem \ref{eq1}, we deduce that
$$\kappa(X, K_X+\Delta+f^*L)=\kappa(X, lD+f^*L)$$
$$=\kappa(F, K_F+\Delta|_F)+\kappa(A, \widehat{\det} (f_*\mathcal{O}_X(lD)\otimes\mathcal{O}_A(L)))$$
$$=\kappa(F, K_F+\Delta|_F)+\dim V^0(A, f_*\mathcal{O}_X(lD)\otimes\mathcal{O}_A(L)).$$
We have $\widehat{\det} (f_*\mathcal{O}_X(lD)\otimes\mathcal{O}_A(L))\cong\widehat{\det} f_*\mathcal{O}_X(lD)\otimes\mathcal{O}_A(rL)$ where $r=\rank f_*\mathcal{O}_X(lD)$. Since $\kappa(A, \widehat{\det} f_*\mathcal{O}_X(lD))\ge0$ by \cite[Lemma 3.3]{Men22} and $\kappa(A, L)\ge0$, we deduce that
$$\kappa(A, \widehat{\det} (f_*\mathcal{O}_X(lD)\otimes\mathcal{O}_A(L)))=\kappa(A, \widehat{\det} f_*\mathcal{O}_X(lD)\otimes\mathcal{O}_A(rL))$$
$$=\kappa(A, \widehat{\det} f_*\mathcal{O}_X(lD)\otimes\mathcal{O}_A(L)).$$
We finish the proof.
\end{proof}

Next, we prove a general theorem and then prove Theorem \ref{eqconc} as an application of the theorem.

\begin{thm}\label{eqcon}
Let $f\colon X \to A$ be a surjective morphism from a klt pair $(X, \Delta)$ to an abelian variety $A$, $m\ge1$ a rational number, and $D$ a Cartier divisor on $X$ such that $D\sim_{\Q}m(K_X+\Delta)$. Let
$$f_*\mathcal{O}_X(D)\cong \bigoplus_{i\in I}(\alpha_i\otimes p_i^*\mathcal{F}_i)$$
be any finite direct sum decomposition of $f_*\mathcal{O}_X(D)$ such that each $A_i$ is an abelian variety, each $p_i\colon A\to A_i$ is a fibration, each $\mathcal{F}_i$ is a nonzero M-regular coherent sheaf on $A_i$, and each $\alpha_i\in\Pic^0(A)$ is a torsion line bundle. If $\kappa(A, \widehat{\det} f_*\mathcal{O}_X(D))=\dim A_j$ for some $j\in I$, then there exists a fibration $t_{ji}\colon A_j\to A_i$ such that $p_i=t_{ji}\circ p_j$ for every $i\in I$.
\[
	\begin{tikzcd}
			(X, \Delta) \dar{f} \\
			A \dar{p_j} \rar{p_i} & A_i \\
			A_j \urar[swap]{t_{ji}}
	\end{tikzcd}
\]
\end{thm}

\begin{proof}
If $f_*\mathcal{O}_X(D)=0$, then the statement is trivial. Thus we can assume $f_*\mathcal{O}_X(D)\neq0$. Let
$$f_*\mathcal{O}_X(D)\cong \bigoplus_{i\in I}(\alpha_i\otimes p_i^*\mathcal{F}_i)$$
be a finite direct sum decomposition of $f_*\mathcal{O}_X(D)$ such that each $A_i$ is an abelian variety, each $p_i\colon A\to A_i$ is a fibration, each $\mathcal{F}_i$ is a nonzero M-regular coherent sheaf on $A_i$, each $\alpha_i\in\Pic^0(A)$ is a torsion line bundle, and
$$\kappa(A, \widehat{\det} f_*\mathcal{O}_X(D))=\dim A_j$$
for some $j\in I$. We consider the following commutative diagram where $b\in A_j$ is a general closed point, and $f_{b}$ is obtained by base change from $f$ via the closed immersion $i_b$ of the closed fiber $A_b$ of $p_j$ over $b$.
\[
	\begin{tikzcd}
	(X_b, \Delta|_{X_b})\dar{f_b} \rar	&	(X, \Delta) \dar{f} \\
	A_b\dar \rar{i_b}	&	A \dar{p_j}\rar{p_i} & A_i\\
	\{b\} \rar	&	A_j
	\end{tikzcd}
\]
We can choose $b$ sufficiently general such that $(X_b, \Delta|_{X_b})$ is klt, and $f_b$ is surjective. By \cite[Lemma 3.2]{Men22}, we can choose $b$ sufficiently general such that
$$i_b^*\widehat{\det}f_*\mathcal{O}_{X}(D)\cong \widehat{\det} {f_b}_*\mathcal{O}_{X_b}(D|_{X_b}).$$
By the same argument as in the proof of Theorem \ref{m2}, we have
$$\widehat{\det} f_*\mathcal{O}_X(D)\cong(\bigotimes_{i\in I}p_i^*\widehat{\det}\mathcal{F}_i)\otimes(\bigotimes_{i\in I}\alpha_i^{\otimes\rank\mathcal{F}_i}),$$
and $\widehat{\det}\mathcal{F}_i$ is an ample line bundle. We have $\kappa(A, p_i^*\widehat{\det}\mathcal{F}_i)\ge0$ for every $i\in I$. Thus we have an injective morphism
$$p_j^*(\widehat{\det}\mathcal{F}_j)^{\otimes N}\hookrightarrow (\widehat{\det} f_*\mathcal{O}_X(D))^{\otimes N}$$
for sufficiently divisible positive integer $N$. By \cite[Proposition 1.14]{Mor87}, we deduce that
$$\kappa(A, \widehat{\det} f_*\mathcal{O}_X(D))=\kappa(A_b, i_b^*\widehat{\det} f_*\mathcal{O}_X(D))+\dim A_j$$
$$=\kappa(A_b, \widehat{\det} {f_b}_*\mathcal{O}_{X_b}(D|_{X_b}))+\dim A_j$$
for very general $b$. Let $T$ be a divisor on $A$ such that $\mathcal{O}_A(T)\cong\widehat{\det} f_*\mathcal{O}_X(D)$. By \cite[Lemma 3.3]{Men22}, we have
$$\kappa(A, \mathcal{O}_A(T))=\kappa(A, \widehat{\det} f_*\mathcal{O}_X(D))\ge0.$$
Thus there exists an effective $\Q$-divisor $E$ on $A$ such that $T\sim_{\Q}E$. We can choose a rational number $\varepsilon>0$ small enough such that $(A, \varepsilon E)$ is a klt pair. By \cite[Theorem 4.2]{HMX18}, $\kappa(A_b, K_{A_b}+\varepsilon E|_{A_b})$ is constant for general $b$. We deduce that
$$\kappa(A_b, \widehat{\det} {f_b}_*\mathcal{O}_{X_b}(D|_{X_b}))=\kappa(A_b, i_b^*\widehat{\det} f_*\mathcal{O}_X(D))=\kappa(A_b, i_b^*\mathcal{O}_A(T))$$
$$=\kappa(A_b, i_b^*E)=\kappa(A_b, \varepsilon E|_{A_b})=\kappa(A_b, K_{A_b}+\varepsilon E|_{A_b})$$
is constant for general $b$. Thus
$$\kappa(A, \widehat{\det} f_*\mathcal{O}_X(D))=\kappa(A_b, i_b^*\widehat{\det} f_*\mathcal{O}_X(D))+\dim A_j$$
$$=\kappa(A_b, \widehat{\det} {f_b}_*\mathcal{O}_{X_b}(D|_{X_b}))+\dim A_j$$
for general $b$. Since $\kappa(A, \widehat{\det} f_*\mathcal{O}_X(D))=\dim A_j$, we deduce that
$$\kappa(A_b, i_b^*\widehat{\det} f_*\mathcal{O}_X(D))=\kappa(A_b, \widehat{\det} {f_b}_*\mathcal{O}_{X_b}(D|_{X_b}))=0$$
for general $b$. Thus by the same argument as in the proof of Theorem \ref{m2}, we have
$$\kappa(A_b, i_b^*p_i^*\widehat{\det}\mathcal{F}_i)=0$$
for every $i\in I$.

We claim that the image of $A_b$ under the morphism $p_i\circ i_b$ is a point for sufficiently general closed point $b$. We prove the claim by contradiction. Assume that it is not a point. Then the scheme theoretic image denoted by $(p_i\circ i_b)(A_b)$ of $p_i\circ i_b$ is an abelian variety of positive dimension since $p_i\circ i_b$ is a morphism between abelian varieties. We consider the following commutative diagram where $v$ is the closed immersion from $(p_i\circ i_b)(A_b)$ to $A_i$, and $u$ is the surjective morphism induced by $p_i\circ i_b$.
\[
	\begin{tikzcd}
	A_b\dar{u} \rar{i_b}	&	A \dar{p_i} \\
	(p_i\circ i_b)(A_b)\rar{v} &	A_i
	\end{tikzcd}
\]
Since $v$ is a closed immersion, we deduce that $v^*\widehat{\det}\mathcal{F}_i$ is an ample line bundle on $(p_i\circ i_b)(A_b)$. Since $u$ is surjective, we have
$$0=\kappa(A_b, i_b^*p_i^*\widehat{\det}\mathcal{F}_i)=\kappa(A_b, u^*v^*\widehat{\det}\mathcal{F}_i)$$
$$=\kappa((p_i\circ i_b)(A_b), v^*\widehat{\det}\mathcal{F}_i)=\dim (p_i\circ i_b)(A_b)>0,$$
which is a contradiction. Thus $(p_i\circ i_b)(A_b)$ must be a point, and we finish proving the claim. Since $A_j$ is smooth and $A_i$ is an abelian variety, we deduce that there exists a morphism $t_{ji}\colon A_j\to A_i$ such that $p_i=t_{ji}\circ p_j$ for every $i\in I$ by \cite[Lemma 14]{Kaw81}. Since $p_j$ and $p_i$ are fibrations, $t_{ji}$ is also a fibration.
\end{proof}

\begin{proof}[Proof of Theorem \ref{eqconc}]
If $f_*\mathcal{O}_X(D)=0$, then the statement is trivial. Thus we can assume $f_*\mathcal{O}_X(D)\neq0$. By \cite[Theorem 1.3]{Men21}, $f_*\mathcal{O}_X(D)$ has the Chen--Jiang decomposition
$$f_*\mathcal{O}_X(D)\cong \bigoplus_{i\in I}(\alpha_i\otimes p_i^*\mathcal{F}_i),$$
where each $A_i$ is an abelian variety, each $p_i\colon A\to A_i$ is a fibration, each $\mathcal{F}_i$ is a nonzero M-regular coherent sheaf on $A_i$, and each $\alpha_i\in\Pic^0(A)$ is a torsion line bundle. By \cite[Lemma 3.3]{LPS20}, we deduce that
$$V^0(A, f_*\mathcal{O}_X(D))=\bigcup_{i\in I}\alpha^{-1}_i\otimes p_i^*\Pic^0(A_i)$$
and
$$\dim V^0(A, f_*\mathcal{O}_X(D))=\max_{i\in I} \dim A_i.$$

Assume that
$$\kappa(A, \widehat{\det} f_*\mathcal{O}_X(D))=\dim V^0(A, f_*\mathcal{O}_X(D)).$$
Thus there exists some $j\in I$ such that
$$\dim A_j=\max_{i\in I} \dim A_i=\dim V^0(A, f_*\mathcal{O}_X(D))=\kappa(A, \widehat{\det} f_*\mathcal{O}_X(D)).$$
Thus by Theorem \ref{eqcon}, there exists a fibration $t_{ji}\colon A_j\to A_i$ such that $p_i=t_{ji}\circ p_j$ for every $i\in I$. We can choose a sufficiently divisible positive integer $M$ such that $\alpha_i^{\otimes M}\cong\mathcal{O}_A$ for every $i\in I$ since each $\alpha_i\in\Pic^0(A)$ is a torsion line bundle. We deduce that
$$V^0(A, f_*\mathcal{O}_X(D))^{\otimes M}=\bigcup_{i\in I}(\alpha^{-1}_i)^{\otimes M}\otimes p_i^*\Pic^0(A_i)^{\otimes M}=\bigcup_{i\in I}p_i^*\Pic^0(A_i)$$
$$=\bigcup_{i\in I}p_j^*t_{ji}^*\Pic^0(A_i)=p_j^*(\bigcup_{i\in I}t_{ji}^*\Pic^0(A_i))=p_j^*\Pic^0(A_j).$$
In particular, $V^0(A, f_*\mathcal{O}_X(D))^{\otimes M}$ is irreducible for sufficiently divisible positive integer $M$.

Assume that $V^0(A, f_*\mathcal{O}_X(D))^{\otimes M}$ is irreducible for sufficiently divisible positive integer $M$. We can choose a sufficiently divisible positive integer $M$ such that $\alpha_i^{\otimes M}\cong\mathcal{O}_A$ for every $i\in I$ and $V^0(A, f_*\mathcal{O}_X(D))^{\otimes M}$ is irreducible. We have that
$$V^0(A, f_*\mathcal{O}_X(D))^{\otimes M}=\bigcup_{i\in I}(\alpha^{-1}_i)^{\otimes M}\otimes p_i^*\Pic^0(A_i)^{\otimes M}=\bigcup_{i\in I}p_i^*\Pic^0(A_i).$$
We deduce that there exists some $j\in I$ such that
$$p_i^*\Pic^0(A_i)\subseteq p_j^*\Pic^0(A_j)$$
for every $i\in I$ since $V^0(A, f_*\mathcal{O}_X(D))^{\otimes M}$ is irreducible. Thus we deduce that there exists a morphism $t_{ji}\colon A_j\to A_i$ such that $p_i=t_{ji}\circ p_j$ for every $i\in I$. Since $p_j$ and $p_i$ are fibrations, $t_{ji}$ is also a fibration.
\[
	\begin{tikzcd}
		A \arrow[r, "p_j"] \arrow[rr, bend right, "p_i"]& A_j \arrow[r, "t_{ji}"] & A_i
	\end{tikzcd}
\]
We deduce that
$$\dim V^0(A, f_*\mathcal{O}_X(D))=\max_{i\in I} \dim A_i=\dim A_j.$$
By the same argument as in the proof of Theorem \ref{m2}, we have
$$\widehat{\det} f_*\mathcal{O}_X(D)\cong(\bigotimes_{i\in I}p_i^*\widehat{\det}\mathcal{F}_i)\otimes(\bigotimes_{i\in I}\alpha_i^{\otimes\rank\mathcal{F}_i}),$$
and $\widehat{\det}\mathcal{F}_i$ is an ample line bundle. We deduce that $t_{ji}^*\widehat{\det}\mathcal{F}_i$ is a nef line bundle, and $\bigotimes_{i\in I}t_{ji}^*\widehat{\det}\mathcal{F}_i$ is an ample line bundle. Thus
$$\kappa(A, \widehat{\det} f_*\mathcal{O}_X(D))=\kappa(A, (\bigotimes_{i\in I}p_i^*\widehat{\det}\mathcal{F}_i)\otimes(\bigotimes_{i\in I}\alpha_i^{\otimes\rank\mathcal{F}_i}))$$
$$=\kappa(A, \bigotimes_{i\in I}p_i^*\widehat{\det}\mathcal{F}_i)=\kappa(A, \bigotimes_{i\in I}p_j^*t_{ji}^*\widehat{\det}\mathcal{F}_i)=\kappa(A_j, \bigotimes_{i\in I}t_{ji}^*\widehat{\det}\mathcal{F}_i)$$
$$=\dim A_j=\dim V^0(A, f_*\mathcal{O}_X(D)).$$
\end{proof}

\section{Results on Conjecture \ref{sic}}\label{4}

In this section, we mainly study Conjecture \ref{sic} and prove Theorems \ref{sicpr} and \ref{eqmad}. We start with a useful lemma.

\begin{lemma}\label{nzmad}
Let $f\colon X \to Y$ be a surjective morphism from a klt pair $(X, \Delta)$ to a smooth projective variety $Y$ of maximal Albanese dimension, $m\geq1$ a rational number, and $D$ a Cartier divisor on $X$ such that $D\sim_{\Q}m(K_X+\Delta)$. If $f_*\mathcal{O}_X(D)\neq0$, then
$$\kappa(Y, \widehat{\det} f_*\mathcal{O}_X(D))\ge0.$$
\end{lemma}

\begin{proof}
Since $Y$ is of maximal Albanese dimension, the Albanese morphism $a_Y\colon Y\to\Alb(Y)$ of $Y$ is generically finite onto its image. By \cite[Theorem 1.1]{Men21}, there exists an isogeny $\varphi\colon A'\to \Alb(Y)$ such that $\varphi^*(a_Y)_*f_*\mathcal{O}_X(D)$ is globally generated where $A'$ is an abelian variety. We consider the following base change diagram.
\[
	\begin{tikzcd}
	X' \dar{f'} \rar	&	X \dar{f} \\
	Y' \dar {a_Y'} \rar{\varphi'}	&	Y \dar{a_Y}  \\
	A' \rar	{\varphi} &	\Alb(Y)
	\end{tikzcd}
\]
Since $a_Y$ is generically finite onto its image, the morphism $a_Y'$ which is obtained by base change from $a_Y$ via $\varphi$ is also generically finite onto its image. Thus we can take a Stein factorization of $a_Y'$ such that $a_Y'=h\circ g$, where $g\colon Y'\to Z$ is birational, $h\colon Z\to A'$ is finite onto its image, and $Z$ is a normal projective variety.
\[
	\begin{tikzcd}
		Y' \arrow[r, "g"] \arrow[rr, bend right, "a_Y'"]& Z \arrow[r, "h"] & A'
	\end{tikzcd}
\]
Since $\varphi^*(a_Y)_*f_*\mathcal{O}_X(D)$ is globally generated, we deduce that 
$$h_*g_*(\varphi')^*f_*\mathcal{O}_X(D)\cong(a_Y')_*(\varphi')^*f_*\mathcal{O}_X(D)\cong\varphi^*(a_Y)_*f_*\mathcal{O}_X(D)$$
is globally generated. We deduce that $h^*h_*g_*(\varphi')^*f_*\mathcal{O}_X(D)$ is globally generated. Since $h$ is a finite morphism, we deduce that the following adjoint morphism
$$h^*h_*g_*(\varphi')^*f_*\mathcal{O}_X(D)\to g_*(\varphi')^*f_*\mathcal{O}_X(D)$$
is surjective and thus $g_*(\varphi')^*f_*\mathcal{O}_X(D)$ is globally generated. We deduce that $g^*g_*(\varphi')^*f_*\mathcal{O}_X(D)$ is globally generated. Since $g$ is a birational morphism, we deduce that the following adjoint morphism
$$g^*g_*(\varphi')^*f_*\mathcal{O}_X(D)\to (\varphi')^*f_*\mathcal{O}_X(D)$$
is generically surjective and thus $(\varphi')^*f_*\mathcal{O}_X(D)$ is generically globally generated. Thus $\widehat{\det}(\varphi')^*f_*\mathcal{O}_X(D)$ is generically globally generated. In particular, $\widehat{\det}(\varphi')^*f_*\mathcal{O}_X(D)$ has nonzero sections. By \cite[Lemma 3.1]{Men22}, we deduce that
$$\kappa(Y, \widehat{\det} f_*\mathcal{O}_X(D))=\kappa(Y', \varphi'^*\widehat{\det} f_*\mathcal{O}_X(D))=\kappa(Y', \widehat{\det}\varphi'^*f_*\mathcal{O}_X(D))\ge0.$$
\end{proof}

First, we deal with a special case of Theorem \ref{eqmad} when the base $Y$ admits a birational morphism to an abelian variety.

\begin{prop}\label{eqbta}
Let $f\colon X \to Y$ be a surjective morphism from a klt pair $(X, \Delta)$ to a smooth projective variety $Y$ which admits a birational morphism $g\colon Y \to A$ to an abelian variety $A$, $F$ the general fiber of $f$, $m>1$ a rational number, and $D$ a Cartier divisor on $X$ such that $D\sim_{\Q}m(K_X+\Delta)$ and $f_*\mathcal{O}_X(D)\neq0$. Let $L$ be a divisor on $Y$ such that $\kappa(Y, L)\ge0$. For every positive integer $l$, we have
$$\kappa(X, K_X+\Delta+f^*L)=\kappa(F, K_F+\Delta|_F)+\kappa(Y, \widehat{\det} f_*\mathcal{O}_X(lD)\otimes\mathcal{O}_Y(K_Y+L)).$$
\end{prop}

\begin{proof}
We only need to prove the equality when $l=1$. Since $f_*\mathcal{O}_X(D)\neq0$, $\kappa(F, D|_F)\ge0$ and thus $\kappa(F, K_F+\Delta|_F)\ge0$. Then $\kappa(X, K_X+\Delta)\ge0$ by Theorem \ref{eq1}. Since $\kappa(Y, L)\ge0$, there exists an effective $\Q$-divisor $L'$ on $Y$ such that $L\sim_{\Q}L'$. We can choose a positive integer $k$ which is sufficiently big such that $(X, \Delta+\frac{f^*L'}{km})$ is klt. We have $kD+f^*L\sim_{\Q}km(K_X+\Delta)+f^*L'\sim_{\Q}km(K_X+\Delta+\frac{f^*L'}{km})$. By Theorem \ref{eq1}, we deduce that
$$\kappa(X, K_X+\Delta+\frac{f^*L'}{km})=\kappa(F, K_F+\Delta|_F)+\kappa(A, \widehat{\det} g_*f_*\mathcal{O}_X(kD+f^*L)).$$
Since $\kappa(X, K_X+\Delta)\ge0$ and $\kappa(X, f^*L)\ge0$, we deduce that $\kappa(X, K_X+\Delta+f^*L)=\kappa(X, K_X+\Delta+\frac{f^*L'}{km})$. Thus
$$\kappa(X, K_X+\Delta+f^*L)=\kappa(F, K_F+\Delta|_F)+\kappa(A, \widehat{\det} g_*f_*\mathcal{O}_X(kD+f^*L)).$$
Since $g$ is a birational morphism, we deduce that
$$\widehat{\det} g_*f_*\mathcal{O}_X(kD+f^*L)\cong (g_*\widehat{\det} f_*\mathcal{O}_X(kD+f^*L))^{\vee\vee}$$
$$\cong(g_*\widehat{\det}(f_*\mathcal{O}_X(kD)\otimes\mathcal{O}_Y(L)))^{\vee\vee}\cong(g_*(\widehat{\det}f_*\mathcal{O}_X(kD)\otimes\mathcal{O}_Y(rL)))^{\vee\vee}$$
$$\cong(g_*\widehat{\det}f_*\mathcal{O}_X(kD))^{\vee\vee}\otimes(g_*\mathcal{O}_Y(rL))^{\vee\vee}\cong\widehat{\det} g_*f_*\mathcal{O}_X(kD)\otimes\mathcal{O}_A(rg_*L),$$
where $r=\rank f_*\mathcal{O}_X(kD)$ and $g_*L$ is a divisor on $A$ which is the pushforward of $L$ under the birational morphism $g$. Since $g_*f_*\mathcal{O}_X(kD)\neq0$, we deduce that $\kappa(A, \widehat{\det} g_*f_*\mathcal{O}_X(kD))\ge0$ by Lemma \ref{nzmad}. Since $\kappa(Y, L)\ge0$, $\kappa(A, g_*L)\ge0$. Thus we deduce that
$$\kappa(A, \widehat{\det} g_*f_*\mathcal{O}_X(kD+f^*L))=\kappa(A, \widehat{\det} g_*f_*\mathcal{O}_X(kD)\otimes\mathcal{O}_A(rg_*L))$$
$$=\kappa(A, \widehat{\det} g_*f_*\mathcal{O}_X(kD)\otimes\mathcal{O}_A(g_*L)).$$
Let $T$ be a divisor on $Y$ such that $\mathcal{O}_Y(T)\cong\widehat{\det}f_*\mathcal{O}_X(D)$. By Theorem \ref{eq2}, we deduce that
$$\kappa(X, K_X+\Delta+f^*g^*g_*L)=\kappa(F, K_F+\Delta|_F)+\kappa(A, \widehat{\det} g_*f_*\mathcal{O}_X(l'D)\otimes\mathcal{O}_A(g_*L))$$
for every positive integer $l'$. In particular, we have
$$\kappa(A, \widehat{\det} g_*f_*\mathcal{O}_X(kD)\otimes\mathcal{O}_A(g_*L))=\kappa(A, \widehat{\det} g_*f_*\mathcal{O}_X(D)\otimes\mathcal{O}_A(g_*L))$$
$$=\kappa(A, (g_*\widehat{\det}f_*\mathcal{O}_X(D))^{\vee\vee}\otimes\mathcal{O}_A(g_*L))=\kappa(A, (g_*\mathcal{O}_Y(T))^{\vee\vee}\otimes\mathcal{O}_A(g_*L))$$
$$=\kappa(A, \mathcal{O}_A(g_*T)\otimes\mathcal{O}_A(g_*L))=\kappa(A, g_*T+g_*L).$$
Since $g$ is a birational morphism to an abelian variety $A$, there exists an effective and $g$-exceptional divisor $E$ on $Y$ such that $K_Y\sim E$ and $\Supp E=\Exc(g)$. There also exist effective and $g$-exceptional divisors $E_1$ and $E_2$ on $Y$ such that $g^*(g_*T+g_*L)+E_1=T+L+E_2$. By Lemma \ref{nzmad}, $\kappa(Y, T)=\kappa(Y, \widehat{\det}f_*\mathcal{O}_X(D))\ge0$. Since we also have $\kappa(Y, L)\ge0$ and $\Supp E=\Exc(g)$, we deduce that
$$\kappa(Y, \widehat{\det} f_*\mathcal{O}_X(D)\otimes\mathcal{O}_Y(K_Y+L))=\kappa(Y, T+E+L)=\kappa(Y, T+E+L+E_2)$$
$$=\kappa(Y, E+g^*(g_*T+g_*L)+E_1)=\kappa(A, g_*T+g_*L)$$
$$=\kappa(A, \widehat{\det} g_*f_*\mathcal{O}_X(kD)\otimes\mathcal{O}_A(g_*L))=\kappa(A, \widehat{\det} g_*f_*\mathcal{O}_X(kD+f^*L)).$$
Thus we have
$$\kappa(X, K_X+\Delta+f^*L)=\kappa(F, K_F+\Delta|_F)+\kappa(A, \widehat{\det} g_*f_*\mathcal{O}_X(kD+f^*L))$$
$$=\kappa(F, K_F+\Delta|_F)+\kappa(Y, \widehat{\det} f_*\mathcal{O}_X(D)\otimes\mathcal{O}_Y(K_Y+L)).$$
\end{proof}

\begin{proof}[Proof of Theorem \ref{eqmad}]
We only need to prove the equality when $l=1$. Since $f_*\mathcal{O}_X(D)\neq0$, $\kappa(F, D|_F)\ge0$ and thus $\kappa(F, K_F+\Delta|_F)\ge0$. Since $Y$ is of maximal Albanese dimension, the Albanese morphism $a_Y\colon Y\to\Alb(Y)$ of $Y$ is generically finite onto its image. Thus we can take a Stein factorization of $a_Y$ such that $a_Y=h\circ g$, where $g\colon Y\to Z$ is birational, $h\colon Z\to\Alb(Y)$ is finite onto its image, and $Z$ is a normal projective variety. By \cite[Theorem 13]{Kaw81}, there exist an \'etale cover $\varphi\colon Z'\to Z$, an abelian variety $A$, and a normal projective variety $W$ such that $Z'\cong A\times W$ and $\kappa(W)=\dim W$. Denote the projection of $Z'$ onto $W$ by $p$. We consider the following commutative diagram, where the morphisms in the middle column are obtained by base change from $g$ and $f$ via $\varphi$ and then $\varphi'$, while on the left we have the respective fibers over a general closed point $w\in W$. We define a $\Q$-divisor $\Delta'$ by $K_{X'}+\Delta'=\varphi''^*(K_X+\Delta)$. Since $\varphi''$ is an \'etale morphism, the new pair $(X', \Delta')$ is klt and $\Delta'$ is effective. Define $D'$ by $\varphi''^*D$ then we have $D'\sim_{\Q}m(K_{X'}+\Delta')$. By the flat base change theorem, we have that $f'_*\mathcal{O}_{X'}(D')\cong\varphi'^*f_*\mathcal{O}_X(D)$.
\[
	\begin{tikzcd}
	(X'_w, \Delta'|_{X'_w})\rar \dar{f'_w} & (X', \Delta') \dar{f'} \rar{\varphi''} & (X, \Delta) \dar{f}  \\
	Y'_w\rar{i_w} \dar{g'_w} & Y' \dar{g'} \rar{\varphi'} & Y \dar{g}\arrow[dd, "a_Y", bend left=40]\\
 	A\cong A_w \dar \rar{} & Z' \dar{p} \rar{\varphi} & Z\dar{h}\\
	\{w\}\rar& W & \Alb(Y)
	\end{tikzcd}
\]
The projective variety $Y'$ is smooth, since $\varphi'$ is \'etale. We can choose $w$ sufficiently general such that $(X'_w, \Delta'|_{X'_w})$ is klt, $Y'_w$ is smooth, $f'_w$ is a fibration, and $g'_w$ is birational. By \cite[Theorem 1.7]{Fuj17} and Theorem \ref{eq1}, we deduce that
$$\kappa(X, K_X+\Delta)=\kappa(X', K_{X'}+\Delta')=\kappa(X'_w, K_{X'_w}+\Delta'|_{X'_w})+\dim W$$
$$\ge\kappa(F, K_F+\Delta|_F)+\dim W\ge0$$
since $\kappa(W)=\dim W$. Since $\kappa(Y, L)\ge0$, there exists an effective $\Q$-divisor $L'$ on $Y$ such that $L\sim_{\Q}L'$. Let $T$ be a divisor on $Y$ such that $\mathcal{O}_Y(T)\cong\widehat{\det}f_*\mathcal{O}_X(D)$. By Lemma \ref{nzmad}, $\kappa(Y, T)=\kappa(Y, \widehat{\det} f_*\mathcal{O}_X(D))\ge0$. Since $\kappa(Y, T)\ge0$, there exists an effective $\Q$-divisor $T'$ on $Y$ such that $T\sim_{\Q}T'$. We can choose a rational number $\varepsilon>0$ sufficiently small such that $(Y', \varepsilon\varphi'^*L'+\varepsilon\varphi'^*T')$ is klt. By \cite[Theorem 1.7]{Fuj17},
$$\kappa(Y, K_Y+\varepsilon L'+\varepsilon T')=\kappa(Y', K_{Y'}+\varepsilon\varphi'^*L'+\varepsilon\varphi'^*T')$$
$$=\kappa(Y'_w, K_{Y'_w}+\varepsilon i_w^*\varphi'^*L'+\varepsilon i_w^*\varphi'^*T')+\dim W$$
since $\kappa(W)=\dim W$. Since $\kappa(Y, K_Y)\ge0$, $\kappa(Y, L)\ge0$, and $\kappa(Y, T)\ge0$, we have $\kappa(Y, K_Y+L+T)=\kappa(Y, K_Y+\varepsilon L'+\varepsilon T')$. Similarly, we have $\kappa(Y'_w, K_{Y'_w}+i_w^*\varphi'^*L+i_w^*\varphi'^*T)=\kappa(Y'_w, K_{Y'_w}+\varepsilon i_w^*\varphi'^*L'+\varepsilon i_w^*\varphi'^*T')$. Thus we have
$$\kappa(Y, K_Y+L+T)=\kappa(Y'_w, K_{Y'_w}+i_w^*\varphi'^*L+i_w^*\varphi'^*T)+\dim W.$$
By \cite[Lemma 3.2]{Men22}, we can choose $w$ sufficiently general such that
$$\mathcal{O}_{Y'_w}(i_w^*\varphi'^*T)\cong i_w^*\varphi'^*\mathcal{O}_{Y}(T)\cong i_w^*\varphi'^*\widehat{\det}f_*\mathcal{O}_X(D)\cong i_w^*\widehat{\det}\varphi'^*f_*\mathcal{O}_X(D)$$
$$\cong i_w^*\widehat{\det}f'_*\mathcal{O}_{X'}(D')\cong\widehat{\det}{f'_w}_*\mathcal{O}_{X'_w}(D'|_{X'_w}).$$
Thus we deduce that
$$\kappa(Y, \widehat{\det} f_*\mathcal{O}_X(D)\otimes\mathcal{O}_Y(K_Y+L))=\kappa(Y, K_Y+L+T)$$
$$=\kappa(Y'_w, K_{Y'_w}+i_w^*\varphi'^*L+i_w^*\varphi'^*T)+\dim W$$
$$=\kappa(Y'_w, \widehat{\det}{f'_w}_*\mathcal{O}_{X'_w}(D'|_{X'_w})\otimes\mathcal{O}_{Y'_w}(K_{Y'_w}+i_w^*\varphi'^*L))+\dim W.$$
By the same argument as above, we deduce that
$$\kappa(X, K_X+\Delta+f^*L)=\kappa(X', K_{X'}+\Delta'+\varphi''^*f^*L)$$
$$=\kappa(X'_w, K_{X'_w}+\Delta'|_{X'_w}+(\varphi''^*f^*L)|_{X'_w})+\dim W$$
$$=\kappa(X'_w, K_{X'_w}+\Delta'|_{X'_w}+f'^*_wi_w^*\varphi'^*L)+\dim W.$$
We have that $D'|_{X'_w}\sim_{\Q}m(K_{X'_w}+\Delta'|_{X'_w})$ and $g'_w$ is birational. By Proposition \ref{eqbta}, we deduce that
$$\kappa(X, K_X+\Delta+f^*L)=\kappa(X'_w, K_{X'_w}+\Delta'|_{X'_w}+f'^*_wi_w^*\varphi'^*L)+\dim W$$
$$=\kappa(F, K_F+\Delta|_F)+\kappa(Y'_w, \widehat{\det}{f'_w}_*\mathcal{O}_{X'_w}(D'|_{X'_w})\otimes\mathcal{O}_{Y'_w}(K_{Y'_w}+i_w^*\varphi'^*L))+\dim W$$
$$=\kappa(F, K_F+\Delta|_F)+\kappa(Y, \widehat{\det} f_*\mathcal{O}_X(D)\otimes\mathcal{O}_Y(K_Y+L)).$$
\end{proof}

\begin{proof}[Proof of Theorem \ref{sicpr}]
Statement (1) follows directly from Theorem \ref{eqmad}.

For statement (2), since $Y$ is of general type, we have that
$$\kappa(X)=\kappa(F)+\dim Y$$
by \cite[Corollary IV]{Vie83a}. By \cite[Theorem III]{Vie83a}, we have that $f_*\omega_{X/Y}^{\otimes m}$ is weakly positive for every positive integer $m$ such that $f_*\omega_{X/Y}^{\otimes m}\neq0$. We deduce that $\widehat{\det} f_* \omega_{X/Y}^{\otimes m}$ is weakly positive and thus pseudo-effective. Since $\omega_Y$ is big, $\widehat{\det} f_* \omega_{X/Y}^{\otimes m}\otimes\omega_Y^{\otimes(mr+1)}$ is also big, where $r=\rank f_* \omega_{X/Y}^{\otimes m}$. We deduce that
$$\kappa(X)=\kappa(F)+\dim Y=\kappa(F)+\kappa(Y, \widehat{\det} f_* \omega_{X/Y}^{\otimes m}\otimes\omega_Y^{\otimes(mr+1)})$$
$$=\kappa(F)+\kappa(Y, \widehat{\det} f_* \omega_{X}^{\otimes m}\otimes\omega_Y).$$

Statement (3) follows from statement (1) since a smooth projective curve with nonnegative Kodaira dimension is of maximal Albanese dimension.
\end{proof}

Next, we prove Theorem \ref{ineqmad}. We start with a special case of Theorem \ref{ineqmad} when the base $Y$ admits a birational morphism to an abelian variety.

\begin{prop}\label{ineqbta}
Let $f\colon X \to Y$ be a surjective morphism from a smooth projective variety $X$ to a smooth projective variety $Y$ which admits a birational morphism $g\colon Y \to A$ to an abelian variety $A$ where $f$ is smooth over an open set $V\subseteq Y$, and $m$ a positive integer. Then
$$\kappa(V)\ge\kappa(Y, \widehat{\det} f_* \omega_{X}^{\otimes m}\otimes\omega_Y).$$
\end{prop}

\begin{proof}
If $f_* \omega_{X}^{\otimes m}=0$, then the statement is trivial. Thus we can assume that $f_* \omega_{X}^{\otimes m}\neq0$. Then $g_*f_* \omega_{X}^{\otimes m}\neq0$. Let $S:=Y\setminus V$ which is a closed subset of $Y$. Then $S=C\cup D$ where $C\subseteq Y$ is a closed set with $\codim_Y C\ge 2$ and $D$ is an effective divisor on $Y$. Since $\kappa(Y)\ge0$, we have that
$$\kappa(V)=\kappa(Y, K_Y+D)$$
by \cite[Lemma 2.6]{MP23}. Since $g$ is a birational morphism to an abelian variety $A$, there exists an effective and $g$-exceptional divisor $E$ on $Y$ such that $K_Y\sim E$ and $\Supp E=\Exc(g)$. The morphism $g$ is an isomorphism after being restricted to the open set $V\setminus\Supp E$. Thus $g\circ f$ is smooth over the open set $g(V\setminus\Supp E)\subseteq A$. By Corollary \ref{ineq1}, we deduce that
$$\kappa(V\setminus\Supp E)=\kappa(g(V\setminus\Supp E))\ge\kappa(A, \widehat{\det} g_*f_* \omega_{X}^{\otimes m}).$$
By \cite[Lemma 2.6]{MP23}, we deduce that
$$\kappa(V\setminus\Supp E)=\kappa(Y, K_Y+(D+E)_{\red})=\kappa(Y, E+(D+E)_{\red})$$
$$=\kappa(Y, D+E)=\kappa(Y, K_Y+D)=\kappa(V).$$
By the same argument as in the proof of Proposition \ref{eqbta}, we have that
$$\kappa(A, \widehat{\det} g_*f_* \omega_{X}^{\otimes m})=\kappa(Y, \widehat{\det} f_* \omega_{X}^{\otimes m}\otimes\omega_Y).$$
Thus we deduce that
$$\kappa(V)=\kappa(V\setminus\Supp E)\ge\kappa(A, \widehat{\det} g_*f_* \omega_{X}^{\otimes m})=\kappa(Y, \widehat{\det} f_* \omega_{X}^{\otimes m}\otimes\omega_Y).$$
\end{proof}
 
\begin{proof}[Proof of Theorem \ref{ineqmad}]
If $f_* \omega_{X}^{\otimes m}=0$, then the statement is trivial. Thus we can assume that $f_* \omega_{X}^{\otimes m}\neq0$. Since $Y$ is of maximal Albanese dimension, the Albanese morphism $a_Y\colon Y\to\Alb(Y)$ of $Y$ is generically finite onto its image. Thus we can take a Stein factorization of $a_Y$ such that $a_Y=h\circ g$, where $g\colon Y\to Z$ is birational, $h\colon Z\to\Alb(Y)$ is finite onto its image, and $Z$ is a normal projective variety. By \cite[Theorem 13]{Kaw81}, there exist an \'etale cover $\varphi\colon Z'\to Z$, an abelian variety $A$, and a normal projective variety $W$ such that $Z'\cong A\times W$ and $\kappa(W)=\dim W$. Denote the projection of $Z'$ onto $W$ by $p$. We consider the following commutative diagram, where the morphisms in the middle column are obtained by base change from $g$ and $f$ via $\varphi$ and then $\varphi'$, while on the left we have the respective fibers over a general closed point $w\in W$. By the flat base change theorem, we have that $f'_*\omega_{X'}^{\otimes m}\cong\varphi'^*f_* \omega_{X}^{\otimes m}$.
\[
	\begin{tikzcd}
	X'_w\rar \dar{f'_w} & X' \dar{f'} \rar{\varphi''} & X \dar{f}  \\
	Y'_w\rar{i_w} \dar{g'_w} & Y' \dar{g'} \rar{\varphi'} & Y \dar{g}\arrow[dd, "a_Y", bend left=40]\\
 	A\cong A_w \dar \rar{} & Z' \dar{p} \rar{\varphi} & Z\dar{h}\\
	\{w\}\rar& W & \Alb(Y)
	\end{tikzcd}
\]
The projective varieties $X'$ and $Y'$ are smooth since $\varphi$ is \'etale. We can choose $w$ sufficiently general such that $X'_w$ and $Y'_w$ are smooth, $f'_w$ is a fibration, and $g'_w$ is birational. Let $S:=Y\setminus V$ which is a closed subset of $Y$. Then $S=C\cup D$ where $C\subseteq Y$ is a closed set with $\codim_Y C\ge 2$ and $D$ is an effective divisor on $Y$. Since $\kappa(Y)\ge0$, we have that
$$\kappa(V)=\kappa(Y, K_Y+D)$$
by \cite[Lemma 2.6]{MP23}. Since $\varphi'$ is finite, we have $\codim_{Y'} \varphi'^{-1}(C)\ge 2$. Thus we can choose $w$ sufficiently general such that $\codim_{Y'_w} i_w^{-1}(\varphi'^{-1}(C))\ge 2$ and $i_w^{-1}(\varphi'^{-1}(D))$ is a divisor. Thus we have
$$\kappa(i_w^{-1}(\varphi'^{-1}(V)))=\kappa(Y'_w, K_{Y'_w}+(i_w^*\varphi'^*D)_{\red})=\kappa(Y'_w, K_{Y'_w}+i_w^*\varphi'^*D)$$
by \cite[Lemma 2.6]{MP23}. We have that $f'_w$ is smooth over $i_w^{-1}(\varphi'^{-1}(V))$. By Proposition \ref{ineqbta}, we deduce that
$$\kappa(i_w^{-1}(\varphi'^{-1}(V)))\ge\kappa(Y'_w, \widehat{\det}{f'_w}_*\omega_{X'_w}^{\otimes m}\otimes\omega_{Y'_w}).$$
By the same argument as in the proof of Theorem \ref{eqmad}, we have that
$$\kappa(V)=\kappa(Y, K_Y+D)=\kappa(Y', K_{Y'}+\varphi'^*D)=\kappa(Y'_w, K_{Y'_w}+i_w^*\varphi'^*D)+\dim W$$
$$=\kappa(i_w^{-1}(\varphi'^{-1}(V)))+\dim W,$$
and
$$\kappa(Y, \widehat{\det} f_* \omega_{X}^{\otimes m}\otimes\omega_Y)=\kappa(Y'_w, \widehat{\det}{f'_w}_*\omega_{X'_w}^{\otimes m}\otimes\omega_{Y'_w})+\dim W$$
for sufficiently general $w$. Thus we deduce that
$$\kappa(V)=\kappa(i_w^{-1}(\varphi'^{-1}(V)))+\dim W\ge\kappa(Y'_w, \widehat{\det}{f'_w}_*\omega_{X'_w}^{\otimes m}\otimes\omega_{Y'_w})+\dim W$$
$$=\kappa(Y, \widehat{\det} f_* \omega_{X}^{\otimes m}\otimes\omega_Y).$$
\end{proof}

\nocite{CP19}
\nocite{PS17}
\nocite{VZ01}

\bibliographystyle{amsalpha}
\bibliography{biblio}

\end{document}